\documentclass{article}
\usepackage{latexsym}
\usepackage[all]{xy}
\usepackage{amssymb,amsmath}

\def\O{\mathcal{O}}

\def\W{\mathcal{W}}

\def\C{\mathcal{C}}
\def\M{\mathcal{M}}

\def\dim{\mathrm{dim}}

\def\separation{\medskip}
  \newtheorem{theorem}{Theorem}[section]

\newcommand{\qed}{\hfill  $\Box$\separation}

\begin{document}

\title{A note on the  Petri loci}
 \author{A. BRUNO -- E. SERNESI\footnote{Both authors are members of GNSAGA-INDAM}}
\date{}
\maketitle \begin{abstract} Let $\M_g$ be the coarse moduli space
of complex projective nonsingular curves of genus $g$. We prove
that when the Brill-Noether number $\rho(g,r,n)$ is non-negative
every component of the Petri locus $P^r_{g,n}\subset \M_g$ whose
general member is a curve $C$ such that $W^{r+1}_n(C) =
\emptyset$, has codimension one in $\M_g$.
\end{abstract}

\section{Introduction}

Let $C$ be a nonsingular irreducible projective curve of genus $g
\geq 2$ defined over $\mathbb{C}$. A pair  $(L,V)$
 consisting of      an invertible sheaf  $L$ on $C$ and of an $(r+1)$-dimensional
 vector  subspace $V \subset H^0(L)$,  $r \ge 0$, is called  \emph{a linear series}
 of dimension $r$  and degree $n=\deg(L)$, or a $g^r_n$. If $V=H^0(L)$
 then the $g^r_n$ is said to be \emph{complete}.

\noindent If $(L,V)$ is a $g^r_n$ then the \emph{Petri map}  for
$(L,V)$ is the natural multiplication map
\[
\xymatrix{ \mu_0(L,V):V \otimes H^0(\omega_C L^{-1}) \ar[r]&
H^0(\omega_C)}
\]
The Petri map for $L$ is
\[
\xymatrix{ \mu_0(L):H^0(L) \otimes H^0(\omega_C L^{-1}) \ar[r]&
H^0(\omega_C)}
\]
Recall that $C$ is called a \emph{Petri curve} if the Petri map
$\mu_0(L)$ is injective for every invertible sheaf $L$ on $C$. By
the Gieseker-Petri theorem \cite{dG82}
 we know that in $\M_g$, the coarse moduli space of nonsingular projective curves of genus $g$,
  the locus of curves which are not Petri is a proper closed subset $P_g$,
  called the \emph{Petri locus}.
  This locus decomposes into several components, according to the numerical types
  and to other properties that linear series can have on a curve of genus $g$.
We will say that  $C$ is \emph{Petri with respect to $g^r_n$'s} if
the Petri map $ \mu_0(L,V)$
   is injective for every $g^r_n$ $(L,V)$ on $C$.

We denote by  $P^r_{g,n}\subset \M_g$ the locus of curves which
are not Petri w.r. to $g^r_n$'s. Then
\[
P_g = \bigcup_{r,n} P^r_{g,n}
\]
where the union is finite by obvious reasons. The structure of
$P^r_{g,n}$ and of $P_g$ is not known in general: both might a
priori have several components and not be equidimensional. In some
special cases  $P^r_{g,n}$ is known to be of pure codimension one
(notably in the obvious case $\rho(g,r,n)=0$,
     and for $r=1$ and $n=g-1$ \cite{mT88}).
     If the
Brill-Noether number
\[
\rho(g,r,n) := g -(r+1)(g-n+r)
\]
is nonnegative then  it is natural to conjecture that $ P^r_{g,n}$
has pure codimension one if it is non-empty.  The evidence is the
fact that $P^r_{g,n}$ is the image in $\M_g$ of a determinantal
scheme $\widetilde{P}^r_{g,n}$ inside the relative Brill-Noether
scheme $\W^r_n\longrightarrow \M_g$, and the expected dimension of
$\widetilde{P}^r_{g,n}$ is $3g-4$.  This is the point of view that
we apply for the proof of our main theorem \ref{T:main1} (see
below). One might ask if even $P_g$ has pure codimension one:
there is not much evidence for this, except that it can be
directly checked to be true for low values of $g$  (see the very recent preprint by M. Lelli-Chiesa \cite{LC11}).

Before stating our result we recall what is known.
 Denote by $\overline{\M}_g$ the moduli space of stable curves, and let
\[
\overline{\M}_g \backslash \M_g =\Delta_0\cup \cdots \cup
\Delta_{[{g\over 2}]}
 \]
 be its boundary, in standard notation.
 In \cite{gF05} G.  Farkas  has   proved the existence of at least one
 divisorial component of $P^1_{g,n}$     in case $\rho(g,1,n) \ge 0$ and $n \le g-1$,
 using the theory of limit linear series.  He found a divisorial component which has a
 nonempty intersection with  $\Delta_1$. Another proof has been given in \cite{CT08},
  by degeneration to a stable curve with $g$ elliptic tails. The
  method of \cite{gF05} has been extended in \cite{gF08} to
  arbitrary $r$.
 In this note without using any degeneration
argument we prove the following result:

\begin{theorem}\label{T:main1}
If $\rho(g,r,n)  \ge 0$ then every component of $P^r_{g,n}$ whose
general member is a curve $C$ such that $W^{r+1}_n(C) =
\emptyset$, has codimension one in $\M_g$.
\end{theorem}

Note that a necessary numerical condition for the existence of a
curve $C$ as in the statement is that $\rho(g,r+1,n) < 0$.  This
condition, together with   $\rho(g,r,n) \ge 0$ gives:
\[
0 \le \rho(g,r,n) < g-n+2(r+1)
\]
or, equivalently:
\[{r\over r+1}g+r \le n < {r+1 \over r+2} g + r+1\]
  For the proof of the theorem we
introduce a modular family $\C \longrightarrow B$ of curves of
genus $g$ (see (i) below for the definition) and we
 use the determinantal description of the relative locus $\W^r_n(\C/B)$ over $B$
 and of the naturally defined closed subscheme $\widetilde{P}^r_{g,n}\subset \W^r_n(\C/B)$
 whose image in $\M_g$ is  $P^r_{g,n}$.  Since it is a determinantal locus, every
 component of $\widetilde{P}^r_{g,n}$ has dimension $\ge 3g-4$. Then a theorem of
 F. Steffen \cite{fS98} ensures that every component of $P^r_{g,n}$ has dimension
 $\ge 3g-4$ as well, thus proving the result.

 In a forthcoming paper (in preparation) we will show the existence of a
 divisorial component of $P^1_{g,n}$ which has a non-empty intersection with $\Delta_0$, when $\rho(g,1,n) \ge 1$.

 \section{Proof of Theorem \ref{T:main1}}

 In this section we fix $g,r,n$ such that $\rho(g,r,n) \ge 0$ and $\rho(g,r+1,n) < 0$.
 Consider the following diagram:
 \begin{equation}\label{E:Jn}
 \xymatrix{
J_n(\C/B)\times_B \C \ar[r] \ar[d]&
 \C \ar[d]^f \\
J_n(\C/B)\ar[r]_-q &B}
 \end{equation}
 where:

 \begin{description}
\item[(i)]   $f$ is a smooth modular family of curves of genus $g$
parametrized by a nonsingular quasi-projective algebraic variety
$B$ of dimension $3g-3$. This means that at each closed point $b
\in B$ the Kodaira-Spencer map $\kappa_b: T_bB \to
H^1(\C(b),T_{\C(b)})$ is an isomorphism. In particular, the
functorial morphism
\[
\xymatrix{ \beta:B \ar[r] & \M_g}
\]
 is finite and dominant. The existence of $f$
is a standard fact, see e.g. \cite{rH10}, Theorem 27.2.

 \item[(ii)]  $J_n(\C/B)$ is the relative Picard variety parametrizing invertible
 sheaves of degree $n$ on the fibres of $f$.

\item[(iii)]  For all closed points $b \in B$ the fibre $\C(b)$
satisfies $W^{r+1}_n(\C(b))=\emptyset$.
  This condition can be satisfied modulo replacing $B$ by an open
  subset
  if necessary, because the condition $W^{r+1}_n(\C(b))=\emptyset$ is open w.r. to $b\in B$.

\item[(iv)] We may even assume that any given specific curve $C$
of genus $g$ satisfying
 $W^{r+1}_n(C)=\emptyset$ appears among the fibres of $f$. In particular we may assume that the dense
 subset Im$(\beta)\subset \M_g$ has a non-empty intersection with all irreducible
 components of $P^r_{g,n}$ whose general element parametrizes a curve $C$ such that  $W^{r+1}_n(C)=\emptyset$.

\end{description}

 Let $\mathcal{P}$ be a Poincar\'e invertible sheaf on  $J_n(\C/B)\times_B \C$.
 Using $\mathcal{P}$ in a well-known way one constructs the relative Brill-Noether scheme
 \[
 \mathcal{W}^r_n(\C/B) \subset J_n(\C/B)
 \]
Consider the restriction of diagram (\ref{E:Jn}) over
$\mathcal{W}^r_n(\C/B)$:
\begin{equation}\label{E:Wrn}
\xymatrix{
 \mathcal{W}^r_n(\C/B)\times_B \C \ar[r]^-{p_2} \ar[d]_{p_1}&
 \C \ar[d]^f \\
\mathcal{W}^r_n(\C/B)\ar[r]_-q &B}
 \end{equation}
 Every irreducible component of $\mathcal{W}^r_n(\C/B)$ has dimension
 $\ge 3g-3+\rho(g,r,n)$ and, since $\rho(g,r,n) \ge 0$, there is a
 component which dominates $B$ \cite{gK72,KL72}. A closed point    $w\in\mathcal{W}^r_n(\C/B)$
 represents  an invertible sheaf $L_w$ on the curve $\C(q(w))$ such that $h^0(L_w) \ge r+1$.
 Denoting again by $\mathcal{P}$ the restriction of $\mathcal{P}$ to
 $\mathcal{W}^r_n(\C/B)\times_B \C$, we have a homomorphism of coherent sheaves on
 $\mathcal{W}^r_n(\C/B)$, induced by multiplication of sections along the fibres of $p_1$:
 \[
 \xymatrix{
 \mu_0(\mathcal{P}):p_{1*}\mathcal{P}\otimes p_{1*}[p_2^*(\omega_{\C/B})\otimes \mathcal{P}^{-1}] \ar[r]&
 p_{1*}[p_2^*\omega_{\C/B}]}
 \]
 By condition (iii) above, these sheaves are locally free, of ranks
 $(r+1)(g-n+r)$ and $g$ respectively. Moreover, by definition, at each point
 $w\in \mathcal{W}^r_n(\C/B)$, the map $\mu_0(\mathcal{P})$ coincides with the Petri map
 \[
 \xymatrix{
 \mu_0(L_w): H^0(\C(q(w)),L_w) \otimes H^0(\C(q(w)),\omega_{\C(q(w))}L_w^{-1})\ar[r]&
 H^0(\C(q(w)),\omega_{\C(q(w))})}
 \]

\emph{Claim:} the vector bundle
\begin{equation}\label{E:vbun}
\left[p_{1*}\mathcal{P}\otimes p_{1*}[p_2^*(\omega_{\C/B})\otimes
\mathcal{P}^{-1}]\right]^\vee \otimes p_{1*}[p_2^*\omega_{\C/B}]
\end{equation}
is $q$-relatively ample.

\emph{Proof of the Claim.} If we restrict diagram (\ref{E:Wrn})
over any $b\in B$ and we let $C=\C(b)$, we obtain:
\[
\xymatrix{ W^r_n(C)\times C \ar[r]^-{\pi_2} \ar[d]_{\pi_1}&
 C  \\
W^r_n(C)}
\]
and the map $\mu_0(\mathcal{P})$ restricts over $W^r_n(C)$ to
\[
\xymatrix{ m_P: \pi_{1*}P \otimes \pi_{1*}[\pi_2^*\omega_C\otimes
P^{-1}] \ar[r]&H^0(C,\omega_C)\otimes \O_{W^r_n}}
\]
where $P = \mathcal{P}_{|W^r_n(C)\times C}$ is a Poincar\'e sheaf
on  $W^r_n(C)\times C$. The dual of the source  of $m_P$ is  an
ample vector bundle (compare \cite{FL81}, \S 2), while the target
is a trivial vector bundle, and therefore
\[
\left[\pi_{1*}P \otimes \pi_{1*}[\pi_2^*\omega_C\otimes
P^{-1}]\right]^\vee \otimes_\mathbb{C}H^0(C,\omega_C)
\]
is an ample vector bundle.  This means that the vector bundle
(\ref{E:vbun}) restricts to an ample vector bundle on the fibres
of $q$.  This implies, by  \cite{aG61}, Th. 4.7.1, applied to the
invertible sheaf $\O(1)$ on the projective bundle associated to
(\ref{E:vbun}), that (\ref{E:vbun}) is $q$-relatively ample. This
proves the Claim.

 Consider the degeneracy scheme:
\[
\widetilde{P}^r_{g,n}:=D_{(r+1)(g-n+r)-1}( \mu_0(\mathcal{P}))
\subset \mathcal{W}^r_n(\C/B)
\]
which is supported on the locus of $w \in \mathcal{W}^r_n(\C/B)$
such that $\mu_0(L_w)$ is not injective. By applying Theorem 0.3
of \cite{fS98} to it we deduce that every irreducible component of
$q(\widetilde{P}^r_{g,n}) \subset B$ has dimension   at least
\[
\dim[\mathcal{W}^r_n(\C/B)] - [g-(r+1)(g-n+r)+1] = 3g-4
\]
 Since $f$ is a modular family, it follows that every
irreducible component of
$\overline{\beta(q(\widetilde{P}^r_{g,n}))} \subset \M_g$ has
dimension $\ge 3g-4$ as well.  But
$\overline{\beta(q(\widetilde{P}^r_{g,n}))} \subset P^r_{g,n} \ne
\M_g$
  and therefore all the components of
$\overline{\beta(q(\widetilde{P}^r_{g,n}))}$ are divisorial.
Since, by (iv), $\overline{\beta(q(\widetilde{P}^r_{g,n}))}$ is
the union of all the components of $P^r_{g,n}$
 whose general element parametrizes a curve $C$ such that $W^{r+1}_n(C)=\emptyset$,
 the theorem   is proved.  \qed

\noindent\textsc{address of the authors:}

\noindent Dipartimento di Matematica,
  Universit\`a Roma Tre \newline Largo S. L. Murialdo 1,
  00146 Roma, Italy.
  \smallskip
\newline   \texttt{bruno@mat.uniroma3.it}
\newline  \texttt{sernesi@mat.uniroma3.it}

\end{document}